\documentclass[11pt]{amsart}
\usepackage[ansinew]{inputenc}
\usepackage{amssymb}
\usepackage{latexsym}
\usepackage{graphics}
\usepackage[dvips]{graphicx}
\usepackage{fancybox}
\usepackage
{amsmath,amsfonts}
\newtheorem{theorem}{Theorem}[section]
\newtheorem{lemma}[theorem]{Lemma}

\newtheorem{remark}{Remark}[section]
\newtheorem{ass}{Assumption}[section]

\begin{document}
\author{ M. M. Cavalcanti }
\address{ Department of Mathematics, State University of
Maring\'a, 87020-900, Maring\'a, PR, Brazil.}
\thanks{Research of Marcelo M. Cavalcanti partially supported by the CNPq Grant
300631/2003-0}

\author{ V. N. Domingos Cavalcanti }
\address{ Department of Mathematics, State University of
Maring\'a, 87020-900, Maring\'a, PR, Brazil.}
\thanks{Research of Valéria N. Domingos Cavalcanti partially supported by the CNPq Grant
304895/2003-2}

\author{ R. Fukuoka }
\address{ Department of Mathematics, State University of
Maring\'a, 87020-900, Maring\'a, PR, Brazil.}

\author{ J. A. Soriano}
\address{ Department of Mathematics, State University of
Maring\'a, 87020-900, Maring\'a, PR, Brazil.}
\thanks{Research of Juan A. Soriano  partially supported by the CNPq Grant
301352/2003-8}

\title[wave on compact manifolds]
{ Uniform Stabilization of the wave equation on compact surfaces
and locally distributed damping}

\maketitle
\begin{abstract}
This paper is concerned with the study of the wave equation on
compact surfaces and locally distributed damping, described by
\begin{equation}
\left.
\begin{array}{l}
u_{tt} -
\Delta_{\mathcal{M}}u+ a(x)\,g(u_{t})=0\, \ \;\,\,\,\,\text{%
on\ \thinspace }\mathcal{M}\times \left] 0,\infty \right[
,\smallskip
\end{array}%
\right.  \nonumber
\end{equation}
where $\mathcal{M}\subset \mathbb{R}^3$ is a smooth (of class $C^3$)
oriented embedded compact surface without boundary, such that
$\mathcal{M}=\mathcal{M}_0 \cup \mathcal{M}_1 $,  where
\begin{eqnarray*}
\mathcal{M}_1:=\{x\in \mathcal{M}; m(x) \cdot \nu(x) > 0\} ~\hbox{
\and }\mathcal{M}_0=\mathcal{M}\backslash \mathcal{M}_1.
\end{eqnarray*}
Here, $m(x):=x-x^0$, ($x^0 \in \mathbb{R}^3$ fixed) and $\nu$ is the
exterior unit normal vector field of $\mathcal{M}.$

\smallskip

For $i=1,\ldots, k$, assume that there exist open subsets
$\mathcal M_{0i} \subset \mathcal M_0$ of $\mathcal M$ with smooth
boundary $\partial \mathcal M_{0i}$ such that $\mathcal M_{0i}$
are umbilical, or more generally, that the principal curvatures
$k_1$ and $k_2$ satisfy $|k_1(x)-k_2(x)|< \varepsilon_i$
($\varepsilon_i$ considered small enough) for all $x\in
\mathcal{M}_{0i}$. Moreover suppose that the {\em mean curvature}
$H$ of each $\mathcal{M}_{0i}$ is {\em non-positive} (i.e. $H\leq
0$ on $\mathcal{M}_{0i}$ for every $i=1,\ldots,k$). If $a(x) \geq
a_0> 0$ on an open subset $\mathcal{M}{\ast} \subset \mathcal M$
that contains $\mathcal{M}\backslash\cup_{i=1}^k \mathcal M_{0i}$
and if $g$ is a monotonic increasing function such that $k |s|
\leq |g(s)| \leq K |s|$ for all $|s| \geq 1$, then uniform decay
rates of the energy hold.
\end{abstract}

\section{Introduction}

\setcounter{equation}{0}

Let $\mathcal{M}$ be a smooth (of class $C^3$) oriented embedded
compact surface without boundary in $\mathbb{R}^{3} $  with
$\mathcal{M} =\mathcal{M} _{0}\cup \mathcal{M}_{1}$, where
\begin{eqnarray}\label{geom. hyp.}
\mathcal{M}_1:=\{x\in \mathcal{M}; m(x) \cdot \nu(x) > 0\} ~\hbox{
\and }\mathcal{M}_0=\mathcal{M}\backslash \mathcal{M}_1.
\end{eqnarray}
Here, $m(x):=x-x^0$, ($x^0 \in \mathbb{R}^3$ fixed) and $\nu$ is the
exterior unit normal vector field of $\mathcal{M}$.

 We denote by $\nabla _{T}$ the tangential-gradient on
$\mathcal{M} $, by $\Delta_{\mathcal{M}}$ the Laplace-Beltrami
operator on $\mathcal{M}$. This paper is devoted to the study of
the uniform stabilization of solutions of the
following damped problem%
\begin{equation}
\left\{
\begin{array}{l}
u_{tt}-\Delta_{\mathcal{M}}u+ a(x)\,g(u_{t})=0\,\,\,\ \ \ \ \ \ \
\ \ \ \ \ \ \ \ \ \ \ \ \ \;\,\,\,\,\text{
on\ \thinspace }\mathcal{M} \times \left] 0,\infty \right[ ,\smallskip \\
u(x,0)=u^0(x),\quad u_t(x,0)=u^1(x)\ \ \ \ \ \ \ \ \ \ \ \ \ \ \ \
x\in \mathcal{M},
\end{array}%
\right.  \label{1.1}
\end{equation}
where $a(x) \geq a_0 > 0$ on an open proper subset
$\mathcal{M}_{\ast}$ of $\mathcal{M}$ and in addition $g$ is a
monotonic increasing function such that $k |s| \leq |g(s)| \leq K
|s|$ for all $|s| \geq 1$.

Stability for the wave equation
\begin{equation}\label{wave-equation}
u_{tt} -\Delta u + f(u) + a(x)\,g(u_t)=0 \hbox{ in }
\Omega\times\mathbb{R}^+,
\end{equation}
where $\Omega$ is a bounded domain in $\mathbb{R}^n$, has been
studied for long time by many authors. When the feedback term
depends on the velocity in a linear way Zuazua \cite{Zuazua}
proved that the energy related to the above equation decays
exponentially if the damping region contains a neighborhood of the
boundary $\partial \Omega$ of $\Omega$ or, at least, contains a
neighborhood $\omega$ of the particular part given by $\{ x\in
\partial \Omega :(x-x_0)\cdot \nu(x)\geq 0\}$. In the same
direction, but when $f=0$,  it is important to mention the work
due to Rauch and Taylor \cite{Rauch} and, subsequently, the
results of Bardos, Lebeau and Rauch \cite{Bardos}, based on
microlocal analysis, that ensures a necessary and sufficient
condition to obtain exponential decay, namely, the damping region
satisfies the well known {\em geometric control condition}. The
classical example of an open subset $\omega$ verifying this
condition  is when $\omega$ is a neighborhood of the boundary.
Later, again considering $f=0$, Nakao \cite{Nakao-1, Nakao-2}
extended the results of Zuazua \cite{Zuazua} treating first the
case of a linear degenerate equation, and then the case of a
nonlinear dissipation $\rho (x,u_t)$ (here, again, $f=0$ was
considered) assuming, as usually, that the function $\rho$ has a
polynomial growth near the origin. Martinez \cite{Martinez}
improved the previous results mentioned above in what concerns the
linear wave equation subject to a nonlinear dissipation $\rho
(x,u_t)$, avoiding the polynomial growth of the function $\rho
(x,s)$ in zero. His proof is based on the piecewise multiplier
technique developed by Liu \cite{Liu} combined with nonlinear
integral inequalities  to show that the energy of the system
decays to zero with a precise decay rate estimate if the damping
region satisfies some geometrical conditions. More recently,  and
still considering $f=0$, Alabau-Boussouira \cite{Alabau} extended
the results due to Martinez \cite{Martinez} by showing {\it
optimal decay rates of energy}. In addition, we would like to
mention the most recent work in this direction due to D. Toundykov
\cite{Daniel} which presents optimal decay rates for solutions to
a semilinear wave equation with localized interior damping and a
source term, subject to Neumann-type boundary condition.

\smallskip

A natural question arises in the context of the wave equation on
compact surfaces: {\em It would be possible to stabilize the
system by considering a localized feedback acting on a part of the
surface ?} In affirmative case, {\em what would be the geometrical
impositions we have to assume on the surface?} When the damping
term acts on the whole surface, the conjecture was studied by
Cavalcanti and Domingos Cavalcanti in \cite{Cavalcanti} and also
by Andrade et al. in \cite{Cavalcanti3, Cavalcanti4} in the
context of viscoelastic problems. However, as far as we are
concerned, there is no result in the literature regarding the {\em
nonlinear wave equation on compact surfaces} when the damping term
acts in a portion $\mathcal{M}_{\ast}$ strictly contained in
$\mathcal{M}$. For the linear case, we can mention the works due
to Rauch\cite{Rauch}, Hitrik \cite{Hitrik} and, more recently
Christianson \cite{Christianson}.

The main goal of this paper is exactly to prove the above
conjecture when the portion of $\mathcal{M}$, where the damping is
effective is strategically chosen. For $i=1,\ldots,k$, assume that
there exist open subsets $\mathcal M_{0i} \subset \mathcal M_0$ of
$\mathcal M$ with smooth boundary $\partial \mathcal M_{0i}$ such
that $\mathcal M_{0i}$ are umbilical, or more generally, that the
principal curvatures $k_1$ and $k_2$ satisfy $|k_1(x)-k_2(x)|<
\varepsilon_i$ ($\varepsilon_i$ considered small enough) for all
$x\in \mathcal{M}_{0i}$. Moreover suppose that the {\em mean
curvature} $H$ of each $\mathcal{M}_{0i}$ is {\em non-positive}
(i.e. $H\leq 0$ on $\mathcal{M}_{0i}$ for every $i=1,\ldots,k$)
and that the damping is effective on an open subset
$\mathcal{M}_{\ast} \subset \mathcal M$ that contains
$\mathcal{M}\backslash\cup_{i=1}^k \mathcal M_{0i}$.

\smallskip

The strategy used to prove the above conjecture is basically to make
use of multipliers and fields as in Lions \cite{Lions1} with new
ingredients that will be clarified in section 4. Indeed, the main
difficulty and the novelty in these kind of problems on surfaces is
how to deal with (or to interpret) the {\em new terms} which appear
in the computations that come from the geometrical structure of
$\mathcal M$. Moreover, this approach can be naturally extended for
semilinear waves where the semilinear function $f(s)$ is assumed to
be super-linear. We would like to emphasize that the proofs of
\cite{Rauch, Bardos,Hitrik}, based on microlocal analysis, {\em do
not extend} to the nonlinear problem (\ref{1.1}). In addition,
making use of arguments due to Lasiecka and Tataru
\cite{Lasiecka-Tataru} we obtain {\it optimal decay rates of the
energy}. The obtained decay rates are optimal, since when we are
able to explicit them (as in Cavalcanti, Domingos Cavalcanti and
Lasiecka  \cite{Cavalcanti5}), they are the same as these optimal
rates derived in the recent works of Alabau-Boussouira \cite{Alabau}
or Toudykov \cite{Daniel}.

Our paper is organized as follows.   Section 2 is concerned with
the statement of the problem  and we introduce some notation . Our
main result is stated in Section 3. Section 4 is devoted to the
proof of the main result.

\section{Statement of Problem}

\setcounter{equation}{0}

Let $\mathcal{M}$ be a smooth oriented embedded compact surface
without boundary in $\mathbb{R}^{3} $ with $\mathcal{M} =\mathcal{M}
_{0}\cup \mathcal{M} _{1}$, where
\begin{eqnarray}\label{2.0}
\mathcal{M}_1:=\{x\in \mathcal{M}; m(x) \cdot \nu(x) > 0\} ~\hbox{
\and }\mathcal{M}_0=\mathcal{M}\backslash \mathcal{M}_1.
\end{eqnarray}
Here, $m$ is the vector field defined by $m(x):=x-x^0$, ($x^0 \in
\mathbb{R}^3$ fixed) and $\nu$ is the exterior unit normal vector
field of $\mathcal{M}$.

In this paper, we investigate the stability properties of
functions $\left[ u(x,t),u_{t}(x,t)\right] $ which solve the
following damped problem:
\begin{equation}
\left\{
\begin{array}{l}
u_{tt}- \Delta_{\mathcal{M}}u+ a(x)\,g(u_{t})=0\, \ \ \ \text{
on\ \thinspace }\mathcal{M}\times \left] 0,\infty \right[ ,\smallskip \\
u(0)=u^0, \quad u_t(0)=u^1,\ \ \ \ \ \
\end{array}%
\right.  \label{3.1}
\end{equation}
where the feedback function $g$ \ satisfies the following
assumptions:

\begin{ass}\label{ass:1}
${}$
\smallskip

$(i) $ $\ \ \ g\left( s\right) $ {\em is continuous and monotone
increasing},

$( ii) $ $\ \ g\left( s\right) s>0$ {\em for} $%
s\neq 0,$

$\left( iii\right) $ $\ \ k\,|s|\leq g\left( s\right) \leq K\,|s
|$ {\em for} $\left\vert s\right\vert >1,$
\end{ass}
\noindent where $k$ and $K$ are two positive constants.

In addition, to obtain the stabilization of problem $(\ref{3.1}),$
we shall need the following geometrical assumption:
\begin{ass}\label{as:2}
Remember that for $i=1,\ldots, k$, $\mathcal M_{0i} \subset \mathcal
M_{0}$ are open sets with smooth boundary $\partial \mathcal M_{0i}$
such that $H\leq 0$ and $\mathcal M_{0i}$ are umbilical
submanifolds, or more generally, that the principal curvatures $k_1$
and $k_2$ satisfy $|k_1(x)-k_2(x)|< \varepsilon_i$ ($\varepsilon_i$
considered small enough) for all $x\in \mathcal{M}_{0i}$. We assume
that $a \in L^{\infty}(\mathcal{M})$ is a nonnegative function such
that
\begin{eqnarray}\label{eq:2.2}
a(x) \geq a_0 > 0, \quad \hbox{\em a. e. on  }\mathcal{M}_{\ast},
\end{eqnarray}
where $\mathcal{M}_{\ast}$ is an open set of $\mathcal{M}$ that
contains $\mathcal{M}\backslash\cup_{i=1}^k \mathcal M_{0i}$.

\end{ass}

In order to fix ideas, Figure 1 shows a compact surface
$\mathcal{M}$ such that there exists only one subset $\mathcal
M_{01}$, which we take as the interior of $\mathcal M_0$.


{\small
\begin{figure}[ht]
\begin{center}
\setlength{\unitlength}{0.73pt}
\begin{picture}(440,210)(0,53)

\put(40,170){$\mathcal{M}_1$} \put(350,155){$\mathcal{M}_{01}$}


\put (460,160){$x_0$} \put(295,234){$\circ$} \put(295,85){$\circ$}
\put(310,237){$\mathcal{M}_{\ast}$}
\put(440,167){\line(-2,1){140}} \put(440,167){\line(2,-1){14}}
\put(455,166){\line(-2,-1){150}} \put(440,165){\line(-5,1){380}}
\put(67,239){\vector(-3,1){12}} \put(319,189){\vector(-4,1){12}}
\put(297,198){$x-x^0$}

\put(344,183){\vector(2,1){25}} \put(350,198){$\nu(x)$}
\put(106,232){\vector(-1,1){25}} \put(20,245){$x-x^0$}
\put(90,255){$\nu(x)$} \put(440,165){\line(5,-1){15}}
\put(445,160){$\bullet$}

\bezier{540}(347,163)(347,252)(207,253)
\bezier{740}(67,163)(67,253)(207,253)
\bezier{740}(347,163)(347,73)(207,73)
\bezier{740}(67,163)(67,73)(207,74)
\put(60,43){\setlength{\unitlength}{0.59pt}
\bezier{740}(9,150)(20,160)(30,165)
\bezier{740}(50,171)(55,172)(70,173)
\bezier{740}(90,174)(100,175)(115,175)
\bezier{740}(140,176)(155,177)(165,177)
\bezier{740}(185,177)(200,177)(210,176)
\bezier{740}(235,175)(250,175)(255,174)
\bezier{740}(280,173)(295,173)(305,172)
\bezier{740}(330,170)(345,168)(345,167)
\bezier{740}(345,167)(350,166)(354,163)
\bezier{740}(9,150)(20,130)(175,130)
\bezier{740}(175,130)(330,130)(354,163)

\bezier{740}(295,56)(300,60)(305,72)
\bezier{740}(312,88)(313,92)(317,102)
\bezier{740}(321,116)(322,120)(324,130)
\bezier{740}(327,144)(328,148)(329,155)

\bezier{740}(329,155)(329,162)(327,169)
\bezier{740}(324,181)(323,184)(315,194)
\bezier{740}(308,210)(306,214)(303,220)
\bezier{740}(298,230)(296,234)(295,238)
\bezier{740}(295,56)(250,150)(295,238)
 }

\put(180,40){Figure 1}
\end{picture}
\end{center}
\caption{The observer is at $x_0$. The subset $\mathcal M_0$ is
the ``visible'' part of $\mathcal M$ and $\mathcal M_1$ is its
complement. The subset $\mathcal{M}_{\ast}\supset \mathcal M-
\cup_{i=1}^k \mathcal M_{0i}=\mathcal M\backslash\mathcal{M}_{01}$
is an open set that contains $\mathcal M\backslash \cup_{i=1}^k
\mathcal M_{0i}$ and the damping is effective there. }

\end{figure}

}

\medskip
In the sequel we define\ by $\Sigma
=\mathcal{M} \times \left] 0,T\right[ ,$ $\Sigma _{i}=\mathcal{M} _{i}\times \left] 0,T%
\right[ ,$ $i=0,1.$

\medskip

Let us considerer the Sobolev spaces $H^s(\mathcal{M})$,\,$s\in
\mathbb{R}$ as in Lions and Magenes \cite{Magenes} section 7.3.

On the other hand, by using the Laplace-Beltrami operator
$\Delta_{\mathcal{M}}$ on $\mathcal{M}$, we can give a more
intrinsic definition of the spaces $H^s(\mathcal{M})$, by
considering
\begin{eqnarray*}
H^{2m}\left( \mathcal{M} \right) &=&\left\{ u\in L^2(\mathcal{M}
)\,/\Delta_{\mathcal{M}}^m \,u \in L^2(\mathcal{M})\right\} ,
\end{eqnarray*}
which, equipped with the canonical norm
\begin{equation}\label{norm H^s II}
\left\Vert u\right\Vert_{H^{2m}(\mathcal{M})}^{2}=\left\Vert
u\right\Vert _{L^{2}(\mathcal{M} )}^{2}+\left\Vert
\Delta_{\mathcal{M}}^m u\right\Vert _{L^{2}(\mathcal{M} )}^{2},
\end{equation}%
is a Hilbert space.
\smallskip

We set
\begin{eqnarray*}
&V:= \{v\in H^1(\mathcal{M}); \int_{\mathcal{M}}
v(x)\,d\mathcal{M} =0 \},&
\end{eqnarray*}
which is a Hilbert space with the topology endowed by
$H^1(\mathcal{M})$.

The condition $\int_{\mathcal{M}} v(x)\,d\mathcal{M} =0$ is
required in order to guarantee the validity of the Poincaré
inequality,
\begin{eqnarray}\label{Poincare}
||f||_{L^2(\mathcal{M})}^2 \leq (\lambda_1)^{-1} ||\nabla_T
f||_{L^2(\mathcal{M})}^2, \quad \hbox{ for all } f\in
H^1(\mathcal{M}),
\end{eqnarray}
where $\lambda_1$ is the first eigenvalue of the Laplace-Beltrami
operator.

We observe that the problem (\ref{3.1}) can be written in the
following form
$$
\frac{dU}{dt} + \mathcal{A}U = G(U_t),
$$
where
\begin{equation*}
\mathcal{A}=\left(
\begin{array}{c}
\,\,0\,\,\,\,\,\,\,\,-I \\
-\Delta_{\mathcal{M}}\,\,\,\,\,\,\,\quad 0
\end{array}
\right)
\end{equation*}
is a maximal monotone operator and $G(\cdot)$ represents a locally
Lipschitz perturbation. So, making use of standard semigroup
arguments we have the following result:

\medskip

\begin{theorem}
${}$
\begin{itemize}
\item $\left( \mathbf{i}\right) $ \textit{Under
the conditions above, problem }$\left( \ref{3.1}\right) $\textit{\ is well posed in%
}$\,\,$\textit{the space}$\,\,V\times L^2(\mathcal{M})$\textit{,\
i.e. for any initial data \thinspace }$\left\{ u^{0},u^{1}\right\}
\in V\times L^2(\mathcal{M})$\textit{, there exists a
unique\thinspace weak solution
of }(\ref{3.1})\textit{\ in the class}%
\begin{equation}
u \in C(\mathbb{R}_{+};V)\cap C^{1}(\mathbb{R}_{+};%
L^2(\mathcal{M})).  \label{3.9}
\end{equation}

\item $ \left(\mathbf{ii}\right) $\textit{In addition, the
velocity term of the solution have the following regularity:}
\begin{equation}
u_{t}\in  L_{loc}^{2}\left( \mathbb{R}_{+};L^{2}\left( \mathcal{M}
\right) \right) , \label{3.10}
\end{equation}%
(\textit{consequently,} $g\left( u_{t}\right) \in
L_{loc}^{2}\left( \mathbb{R}_{+};L^{2}\left( \mathcal{M}\right)
\right) $ by Assumption \ref{ass:1}.
\end{itemize}
\noindent \textit{Furthermore, if} $\left\{ u^{0},u^{1}\right\}
\in \left\{ V\cap H^{2}\left( \mathcal{M} \right) \times V
\right\} $ \textit{then the solution has the following regularity}
\begin{eqnarray*}
 u \in
L^{\infty }\left( \mathbb{R}_{+};V\cap H^{2}\left( \mathcal{M}
\right) \right) \cap W^{1,\infty }\left(\mathbb{R}_{+}; V
\right)\cap W^{2,\infty }\left( \mathbb{R}_{+};L^{2}\left(
\mathcal{M} \right) \right) .
\end{eqnarray*}
\end{theorem}
\smallskip

Supposing that $u $ is the unique global weak solution of problem
(\ref{3.1}), we define the corresponding energy functional by
\begin{equation}
E(t)=\frac{1}{2}\int_{\mathcal{M}}\left[ \left\vert
u_{t}(x,t)\right\vert ^{2}+\left\vert \nabla _T  u(x,t)\right\vert
^{2}\right] d\mathcal{M} . \label{3.11}
\end{equation}

For every solution of (\ref{3.1}) in the class (\ref{3.9}) the
following identity holds%
\begin{equation}
E(t_{2})-E(t_{1})=-\int_{t_{1}}^{t_{2}}\int_{\mathcal{M}
}a(x)\,g(u_{t})u_{t}\,d\mathcal{M} dt, ~\hbox{ for all }
t_{2}>t_{1}\geq 0, \label{3.12}
\end{equation}%
and therefore the energy is a non increasing function of the time
variable $ t $.

\section{Main Result}

\setcounter{equation}{0}

Before stating our stability result, we will define some needed
functions. For this purpose, we are following the ideas firstly
introduced in Lasiecka and Tataru \cite{Lasiecka-Tataru}. For the
reader's comprehension we will repeat them briefly.
Let $h$ be a\ concave, strictly increasing function, with $%
h\left( 0\right) =0$, and such that
\begin{equation}
h\left( s\,g(s))\right) \geq s^{2}+g^{2}(s),\text{\ for
}\left\vert s\right\vert \leq 1.  \label{4.1'}
\end{equation}

Note that such function can be straightforwardly constructed,
given the hypotheses on $g$ in Assumption \ref{ass:1}. With this
function, we define
\begin{equation}
r(.)=h(\frac{.}{meas\left( \Sigma _{1}\right) }).  \label{4.2'}
\end{equation}
As $r$ is monotone increasing, then $cI+r$ is invertible for all
$c\geq 0.$ For $L$ a positive constant, we set
\begin{equation}
p(x)=(cI+r)^{-1}\left( Lx\right) ,  \label{4.3'}
\end{equation}%
where the function $p$ is easily seen to be positive, continuous
and strictly increasing with $p(0)=0$. Finally, let
\begin{equation}
q(x)=x-(I+p)^{-1}\left( x\right) .  \label{4.4'}
\end{equation}%
We can now proceed to state our stability result.

\medskip

\begin{theorem}\label{Theo. 3.1}
Assume that Assumptions \ref{ass:1} and Assumption \ref{as:2} are in
place. Let $u$ be the weak solution of the problem (\ref{3.1}). With
the energy $E(t)$ defined as in (\ref{3.11}), there then exists a
$T_{0}>0$ such that
\begin{equation}
E(t)\leq S\left( \frac{t}{T_{0}}-1\right) ,\text{ \ }\forall
t>T_{0}, \label{4.5'}
\end{equation}%
with $\underset{t\rightarrow \infty }{\lim }S(t)=0,$ where the
contraction semigroup $S(t)$ is the solution of the differential
equation
\begin{equation}
\frac{d}{dt}S(t)+q(S(t))=0,\text{\ \ }S(0)=E(0),  \label{4.6'}
\end{equation}
(where $q$ is as given in (\ref{4.4'})). Here, the constant $L$
(from definition (\ref{4.3'})) will depend on $meas(\Sigma)$, and
the constant $c$(from definition (\ref{4.3'})) is taken here to be $c\equiv \frac{%
k^{-1}+K}{meas\left( \Sigma\right) (1+||a||_{\infty})}.$
\end{theorem}
\smallskip

\noindent \begin{remark}
If the feedback is linear, e. g., $%
g(s)=s,$ then, under the same assumptions as in Theorem \ref{Theo.
3.1}, we have that the energy of problem (\ref{3.1}) decays
exponentially with respect to the initial energy, there exist two
positive constants $C>0$ and $k>0$ such
that%
\begin{equation}
E(t)\leq Ce^{-kt}E(0),\text{ \ \ }t>0.  \label{4.7'}
\end{equation}

\smallskip

As another example, we can consider $g(s) = s^ p $, $ p> 1 $ at
the origin. Since the function $s^{\frac{p+1}{2}} $ is convex for
$p \geq 1 $, then solving
\begin{equation}\label{S1}
S_t + S^{\frac{p+1}{2}} =0,
\end{equation}
we obtain the following polynomial decay rate:
$$ E(t) \leq C(E(0))[E(0)^{\frac{-p +1}{2}} +
t(p-1)]^{\frac{ 2}{-p +1}}.$$

We can find more interesting explicit decay rates in Cavalcanti,
Domingos Cavalcanti and Lasiecka \cite{Cavalcanti5}.
\end{remark}

\section{Proof of Main result}

\setcounter{equation}{0}

\subsection{Preliminaries}
${}$
\smallskip

\smallskip

We collect, below, some few formulas to be invoked in the sequel.

\smallskip

Let $\nu$ be the exterior normal vector field on $\mathcal{M}$. For
all $x\in \mathcal{M}$, we denote by $\pi(x)$ the orthogonal
projection on the tangent plane $T_{x}\mathcal{M}$. Any regular
vector field $q:\mathbb{R}^3 \rightarrow \mathbb{R}^3$ will be split
up as follows:
\begin{equation}
q(x)=q_{T}+(q(x)\cdot\nu(x))\nu(x), \label{4.1}
\end{equation}
where $q_{T}=\pi(x)q(x)$ is the $tangential\, component$ of $q$.

If $\varphi:\mathbb{R}^3 \rightarrow \mathbb{R}$ is a regular
function, we have
\begin{equation}
\nabla\varphi=\partial_{\nu}\varphi\nu+\nabla_{T}\varphi\quad
\hbox{on}\;\;\mathcal{M}, \label{4.2}
\end{equation}
\begin{equation}
 |\nabla\varphi|^2=|\partial_{\nu}\varphi|^2+|\nabla_{T}\varphi|^2 \quad
\hbox{on}\;\;\mathcal{M}, \label{4.3}
\end{equation}%
where $\partial_{\nu}$, is the normal derivative towards the
exterior of $\mathcal{M}$ and $\nabla_{T}\varphi$, is the {\it
tangential gradient} of $\varphi$.

The {\em Laplace- Beltrami operator} $\Delta _{\mathcal{M}}$ of a
function $\varphi :\mathcal{M}\rightarrow \mathbb{R}$ of class
$C^{2}$ is defined by
\begin{equation}
\Delta _{\mathcal{M}}\varphi :=div_{T}\nabla _{T}\varphi ,
\label{4.4}
\end{equation}
where $div_{T}\nabla_{T}\varphi$, is the $divergent$ of the vector
field $\nabla_{T}\varphi$.

Assuming that $\varphi :\mathcal{M}\rightarrow \mathbb{R}$ is a
function of class $C^{1}$ and $q:\mathbb{R}^{3}\rightarrow
\mathbb{R}^{3}$ be a vector field of class $C^{1}$, we have,
\begin{eqnarray}
\int_{\mathcal{M} }q_{T}\cdot \nabla _{T}\varphi \,d\mathcal{M}
&=&-\int_{\mathcal{M} }divq_{T}\,\varphi \,d\mathcal{M},
\label{4.5}
\end{eqnarray}
\begin{eqnarray}
2\varphi (q_T \cdot \nabla_T\varphi)= q_T \cdot
\nabla_T(\varphi^2).\label{E}
\end{eqnarray}

From (4.5) and (4.6), we conclude the following formula
\begin{equation}
2\int_{\mathcal{M} }\varphi (q_{T}\cdot \nabla _{T}\varphi
)\,d\mathcal{M} =\int_{\mathcal{M} }q_{T}\cdot \nabla _{T}(\varphi
^{2})\,d\mathcal{M} =-\int_{\mathcal{M} }div_{T}\,q_{T}|\varphi
|^{2}d\mathcal{M} . \label{second formula}
\end{equation}

We observe that in the particular case when $m(x)=x-x^{0}$, $x\in
\mathbb{R} ^{3}$ and $x^{0}\in \mathbb{R}^{3}$ is a fixed point in
$\mathbb{R}^{3}$, we have
\begin{equation}
\nabla \cdot m = 3,\quad div_{T}\,m_{T}=2+\,(m\cdot \nu )%
Tr\, B .  \label{divergent of mT}
\end{equation}
where $B$ is the second fundamental form of $\mathcal{M}$ (the shape
operator) and $Tr$ is the trace. Let $\varphi $ and $m$ defined as
above. We also have,
\begin{equation}
\nabla _{T}\varphi \cdot \nabla _{T}m_{T}\cdot \nabla _{T}\varphi
=|\nabla _{T}\varphi |^{2}+(m\cdot \nu )(\nabla _{T}\varphi \cdot
B\cdot \nabla _{T}\varphi ).  \label{third formula}
\end{equation}

The proof of the above formulas  can be found in \cite{Nedelec},
\cite{Lemrabet1}, \cite{Heminna3} and references therein.

\smallskip

\begin{remark}\label{remak4.1} The sign of $B$ can change in the
literature. In our case, we remember that $B=-dN$, where $N$ is the
Gauss map with respect to $\nu$.

The formula (\ref{divergent of mT}) can be rewritten  by
\begin{equation}
\nabla \cdot m = 3,\quad div_{T}\,m_{T}=2+2H\,(m\cdot \nu )%
.  \label{divergent of mT3}
\end{equation}
where $H=\frac{trB}{2}$ is the mean curvature of $\mathcal M$.
\end{remark}
\smallskip

We define a continuous linear operator $-\Delta
_{\tilde{\mathcal{M}}}:H^{1}(\tilde{\mathcal{M}})\rightarrow
(H^{1}(\tilde{\mathcal{M}}))^{\prime }$, where
$\tilde{\mathcal{M}}$ is a nomempty open subset of $\mathcal{M}$
(sometimes the whole $\mathcal{M}$) such that

\begin{equation}
\langle -\Delta _{\tilde{\mathcal{M}}}\varphi ,\psi \rangle
=\int_{\tilde{\mathcal{M}}}\nabla _{T}\varphi \cdot \nabla
_{T}\psi \,d\mathcal{M} , \quad \forall \varphi, \psi \in
H^{1}(\tilde{\mathcal{M}}) \label{4.12'}
\end{equation}%
and, in particular,
\begin{equation}
\langle -\Delta _{\tilde{\mathcal{M}}}\varphi ,\varphi \rangle
=\int_{\tilde{\mathcal{M}}}|\nabla _{T}\varphi
|^{2}\,d\mathcal{M}, \quad \forall \varphi \in
H^{1}(\tilde{\mathcal{M}}). \label{4.13'}
\end{equation}

\smallskip

The operator $-\Delta_{\tilde{\mathcal{M}}}+I$ defines an
isomorphism from $H^1(\tilde{\mathcal{M}})$ over
$[H^1(\tilde{\mathcal{M}})]'$. We observe that when
$\tilde{\mathcal{M}}$ is a manifold without boundary, and this is
the case, for instance, if $\tilde{\mathcal{M}}=\mathcal{M}$, we
have $H^1(\tilde{\mathcal{M}})= H_0^1(\tilde{\mathcal{M}})$ and,
consequently,
$[H^1(\tilde{\mathcal{M}})]'=H^{-1}(\tilde{\mathcal{M}})$.

\begin{remark}
It is convenient to observe that all the classical formulas above
stated can be extended for Sobolev spaces by using of density
arguments.
\end{remark}

The proof of Theorem 3.1 proceeds through several steps.

\subsection{An identity}

We begin by proving the following proposition

\smallskip

\noindent \textbf{Proposition 4.2.1. }\textit{Let}
$\mathcal{M}\subset \mathbb{R}^3 $ \textit{be oriented regular
compact surface without boundary} \textit{\ and }$q$ \textit{a
vector field with }$q=q_{T}+(q\cdot\nu) \nu.$\textit{\ Then, for
every regular solution }$u $\textit{\thinspace of
}(\ref{1.1})\textit{ we have the following identity}
\begin{eqnarray}
&&\left[ \int_{\mathcal{M}}u_{t}\, q_{T}\cdot \nabla_T
u\,d\mathcal{M}\right]
_{0}^{T}+\frac{1}{2}\int_{0}^{T}\int_{\mathcal{M}
}(div_{T}q_{T})\left\{ \left\vert u_{t}\right\vert ^{2}-\left\vert
\nabla _{T}
u\right\vert^{2}\right\} d\mathcal{M} dt \label{5.2.1} \\
&&+\int_{0}^{T}\int_{\mathcal{M}} \nabla_T u \cdot \nabla_T q_T
\cdot \nabla_T u \,d\mathcal{M} dt +\int_{0}^{T}\int_{\mathcal{M}
}a(x)\,g(u_{t})(q_{T}\cdot\nabla_T u)d\mathcal{M} dt=0. \nonumber
\end{eqnarray}
\noindent \textbf{Proof. } Multiplying the  equation of
(\ref{1.1}) by the multiplier $q_T \cdot \nabla_T u$ and
integrating on $\mathcal{M}\times ]0,T[$, we obtain
\begin{eqnarray}\label{IDEN}
0= \int_0^T \int_{\mathcal{M}} (u_{tt} - \Delta_{\mathcal{M}} u +
a(x)g(u_t))(q_T \cdot \nabla_T u)\,d\mathcal{M}\,dt.
\end{eqnarray}

Next, we will estimate some terms on the RHS of identity
(\ref{IDEN}).  Taking (\ref{4.12'}), (\ref{E}) and (\ref{second
formula}) into account, we obtain
\begin{eqnarray}
&&\int_{0}^{T}\int_{\mathcal{M}}\left( -\Delta_{\mathcal{M}}
u\right) (q_T\cdot \nabla_T u)\,d\mathcal{M} dt
=\int_{0}^{T}\int_{\mathcal{M}} \nabla_{T} u \cdot\nabla
_{T}(q_T\cdot \nabla_T u)\, d\mathcal{M}
dt   \label{5.2.4} \\
&&=\int_{0}^{T}\int_{\mathcal{M}}\nabla _{T} u \cdot \nabla _{T}
q_T \cdot \nabla_{T} u \,d\mathcal{M} dt
+\frac{1}{2}\int_{0}^{T}\int_{\mathcal{M}}q_T \cdot \nabla_T
[|\nabla_T u|^2]d\mathcal{M}
dt  \notag \\
&&=\int_{0}^{T}\int_{\mathcal{M}}\nabla _{T} u \cdot \nabla _{T}
q_T \cdot \nabla_{T} u \,d\mathcal{M} dt
-\frac{1}{2}\int_{0}^{T}\int_{\mathcal{M}}\left\vert  \nabla _{T}
u\right\vert^{2}div_{T}q_{T}\,d\mathcal{M} dt,  \notag
\end{eqnarray}
and, integrating by parts and considering (\ref{second formula}),
we obtain
\begin{eqnarray}
&&\int_{0}^{T}\int_{\mathcal{M}}\left( u_{tt}+ a(x)\,g( u_{t})
\right)
(q_T \cdot \nabla_T u)\,d\mathcal{M} dt \label{5.2.5}\\
&=&\left[ \int_{\mathcal{M}} u_t(q_T \cdot \nabla_T u)\right]
_{0}^{T}- \int_{0}^{T}\int_{\mathcal{M}}u_t (q_T\cdot \nabla_T u_t)d\mathcal{M} dt  \notag \\
&&+\int_{0}^{T}\int_{\mathcal{M}}a(x)\,g\left(
u_{t}\right) (q_T \cdot \nabla_T u)d\mathcal{M} dt  \notag \\
&=&\left[ \int_{\mathcal{M}}u_{t} (q_T \cdot \nabla_T u)\right]
_{0}^{T}+\frac{1}{2} \int_{0}^{T}\int_{\mathcal{M}
}(div_{T}q_{T})\left\vert u_{t}\right\vert
^{2}d\mathcal{M} dt  \notag \\
&&+\int_{0}^{T}\int_{\mathcal{M}}a(x)\,g\left( u_{t}\right) (q_T
\cdot \nabla_T u)d\mathcal{M} dt. \notag
\end{eqnarray}

Combining (\ref{IDEN}), (\ref{5.2.4}) and ( \ref{5.2.5}), we deduce
(\ref{5.2.1}), which concludes the proof of Proposition 4.2.1.
$\quad \square $

\medskip

Employing (\ref{5.2.1}) with $q(x)= m(x)= x-x^0$ for some $x^0 \in
\mathbb{R}^3$ fixed and taking  (\ref{divergent of mT}) and
(\ref{third formula}) into account, we infer
\begin{eqnarray}
&&\left[ \int_{\mathcal{M}}u_{t}\, m_{T}\cdot \nabla_T
u\,d\mathcal{M}\right] _{0}^{T} +\int_{0}^{T}\int_{\mathcal{M}
}\left\{ \left\vert u_{t}\right\vert ^{2}-\left\vert
\nabla _{T} u\right\vert^{2}\right\} d\mathcal{M} dt \label{first inequality}\\
&&+\int_{0}^{T}\int_{\mathcal{M}} [|\nabla_T u|^2 +(m\cdot \nu)
(\nabla_T u \cdot B \cdot \nabla_T u)]\,d\mathcal{M} dt  \notag \\
&&+\int_{0}^{T}\int_{\mathcal{M}}(m\cdot \nu)H\left\{ \left\vert
u_{t}\right\vert ^{2}-\left\vert \nabla _{T}
u\right\vert^{2}\right\} d\mathcal{M} dt  \notag \\
&&+\int_{0}^{T}\int_{\mathcal{M}}a(x)\,g(u_{t})(m_{T}\cdot\nabla_T
u)d\mathcal{M} dt=0. \nonumber
\end{eqnarray}

\smallskip

\noindent We have the following identity:

\smallskip

\noindent \textbf{Lemma 4.2.3. }\textit{Let }$u $\textit{\ be a
weak solution to problem }(\ref{1.1}) \textit{and} $\xi\in
C^1(\mathcal{M})$.\textit{\ Then}
\begin{eqnarray}\qquad
\left[\int_{\mathcal{M}}u_t\,\xi\, u\,d\mathcal{M}  \right]_0^T
&=& \int_0^T\int_{\mathcal{M}} \xi |u_t|^2 d\mathcal{M} dt-
\int_0^T \int_{\mathcal{M}} \xi |\nabla_T u|^2 d\mathcal{M} dt \label{identity II}\\
&-&  \int_0^T \int_{\mathcal{M}} (\nabla_T u \cdot \nabla_T
\xi)u\,d\mathcal{M} dt -\int_0^T \int_{\mathcal{M}}
a(x)\,g(u_t)\,\xi\,u\,d\mathcal{M} dt.\nonumber
\end{eqnarray}
\textbf{Proof:} Multiplying the first equation of (\ref{1.1}) by
$\xi \,u$ and integrating by parts we obtain the desired. \quad
$\square$
\smallskip

Substituting $\xi=\frac{1}{2}$ in (\ref{identity II}) and
combining the obtained result with identity (\ref{first
inequality}) we deduce
\begin{eqnarray}
&&\left[ \int_{\mathcal{M}}u_{t}\, m_{T}\cdot \nabla_T
u\,d\mathcal{M}\right] _{0}^{T}+\frac{1}{2}\left[
\int_{\mathcal{M}}u_{t}\, u\,d\mathcal{M}\right] _{0}^{T} \label{second inequality}\\
&&+  \int_0^T E(t)\,dt+\int_{0}^{T}\int_{\mathcal{M}
}a(x)\,g(u_{t})(m_{T}\cdot\nabla_T
u)d\mathcal{M} dt\nonumber \\
&& +\frac{1}{2}\int_0^T \int_{\mathcal{M}}a(x)\,g(u_t)u\,
d\mathcal{M} dt\nonumber\\
&&=-\int_{0}^{T}\int_{\mathcal{M}}(m\cdot \nu)H\left\{ \left\vert
u_{t}\right\vert ^{2}-\left\vert
\nabla _{T} u\right\vert^{2}\right\} d\mathcal{M} dt.\notag\\
&&-\int_0^T \int_{\mathcal{M}} (m\cdot \nu) (\nabla_T u \cdot B
\cdot \nabla_T u) \,d\mathcal{M} dt.\nonumber
\end{eqnarray}

\medskip
\noindent{\bf Analysis of the terms which involve the shape operator
$B$}
\medskip

Let us focus our attention on the shape operator $B:T_x \mathcal M
\rightarrow T_x\mathcal M$. There exist an orthonormal basis
$\{e_1,e_2\}$ of $T_x\mathcal M$ such that $Be_1=k_1e_1$ and
$Be_2=k_2e_2$. $k_1$ and $k_2$ are the principal curvatures of
$\mathcal M$ at $x$. The matrix of $B$ with respect to the basis
$\{e_1,e_2\}$ is given by

\begin{equation*}
B:=\left(
\begin{aligned}
&k_1  \quad 0~\\
&~0 \quad k_2
\end{aligned}
\right).
\end{equation*}

Setting $\nabla_T u =(\xi, \eta)$ the coordinates of $\nabla_T u$
in the basis $ \{e_1,e_2\}$, for each $x\in \mathcal{M}$, we
deduce that
\begin{eqnarray}
&&\nabla_T u \cdot B \cdot \nabla_T u = k_1 \xi^2 + k_2
\eta^2.\label{B} \label{trace}
\end{eqnarray}
Then, from (\ref{B}),  we infer
\begin{eqnarray}\label{4.22}
&&(m\cdot \nu) \left[(\nabla_T u \cdot B \cdot \nabla_T u)-\frac12
Tr(B) |\nabla_T u|^2\right]
\\
&&=(m\cdot \nu)\left[ \frac{(k_1-k_2)}{2}\xi^2 +
\frac{(k_2-k_1)}{2}\eta^2\right].\nonumber
\end{eqnarray}

\medskip
\begin{remark}
Observe that this is the precise moment that the intrinsic
properties of the manifold $\mathcal{M}$ appear, that is, we
strongly need that the term $-\int_0^T \int_{\mathcal{M}}(m\cdot
\nu) H u_t^2 \,d\mathcal{M}\,dt$ lies in a region where the
damping term is effective. Remember that the damping term is
effective on an open set $\mathcal{M}_{\ast}$ which contains
$\mathcal M \backslash\cup_{i=1}^k \mathcal M_{0i}$. So, assuming
that $H\leq 0$ and since $m(x) \cdot \nu(x) \leq 0$ on
$\mathcal{M}_0$, we have
\begin{eqnarray*}
-\int_{0}^{T}\int_{\mathcal{M}_0}(m\cdot \nu)H \left\vert
u_{t}\right\vert ^{2} d\mathcal{M} dt\leq 0.
\end{eqnarray*}

In addition, supposing that $\mathcal{M}_{0i}$ is umbilical for
every $i=1,\ldots, k$, then, having (\ref{4.22}) in mind, we also
have that
\begin{eqnarray*}
\int_0^T \int_{\mathcal{M}_{0i}} (m\cdot \nu) \left[H |\nabla_T u|^2
- (\nabla_T u \cdot B \cdot \nabla_T u) \right] \,d\mathcal{M} dt=0,
\; \; i=1,\ldots,k.
\end{eqnarray*}

More generally, assuming that the principal curvatures $k_1$ and
$k_2$ satisfy $|k_1(x)-k_2(x)|< \varepsilon_i$ (here,
$\varepsilon_i$ is assumed sufficiently small) for all $x\in
\mathcal{M}_{0i}$, $i=1, \cdots,k$, we deduce that
\begin{eqnarray*}
&& \left\vert \sum_{i=1}^k \int_0^T \int_{\mathcal{M}_{0i}} (m\cdot
\nu) \left[H |\nabla_T u|^2-(\nabla_T u \cdot B \cdot \nabla_T u)
\right]
\,d\mathcal{M} dt \right\vert \\
&& \leq \sum_{i=1}^k \int_0^T \int_{\mathcal{M}_{0i}} |(m\cdot
\nu)| |k_1-k_2| |\xi^2 + \eta^2| d\mathcal{M}\,dt
\\
&& \leq  \sum_{i=1}^k R_i \varepsilon_i \int_0^T
\int_{\mathcal{M}_{0i}} |\nabla_T u|^2 d\mathcal{M}\,dt \leq  2
\sum_{i=1}^k R_i \varepsilon_i \int_0^T E(t)\,dt,
\end{eqnarray*}
where $R_i = max_{x \in
\overline{\mathcal{M}_{0i}}}||x-x^0||_{\mathbb{R}^3}$.
\end{remark}

Set $\mathcal M_2=\mathcal M\backslash\cup_{i=1}^k \mathcal
M_{0i}$. In the case where $\mathcal M_{0i}$ are umbilical,
recalling (\ref{second inequality}) taking (\ref{4.22}) and Remark
4.3 into consideration, we deduce
\begin{eqnarray}\label{4.24}
\int_0^T E(t)\,dt &\leq& -\left[ \int_{\mathcal{M}}u_{t}\,
m_{T}\cdot \nabla_T u\,d\mathcal{M}\right]
_{0}^{T}-\frac{1}{2}\left[ \int_{\mathcal{M}}u_{t}\,
u\,d\mathcal{M}\right] _{0}^{T}\\
&+&\int_0^T \int_{\mathcal M_2} (m\cdot \nu) \left[H |\nabla_T
u|^2-(\nabla_T u \cdot B \cdot \nabla_T u)\right]
\,d\mathcal{M} dt\nonumber\\
&-&\int_{0}^{T}\int_{\mathcal{M}_2}(m\cdot \nu)H
\left\vert u_{t}\right\vert ^{2} d\mathcal{M} dt\nonumber\\
&-&\int_{0}^{T}\int_{\mathcal{M}
}a(x)\,g(u_{t})(m_{T}\cdot\nabla_T u)d\mathcal{M} dt\nonumber\\
&+&\frac{1}{2}\int_0^T \int_{\mathcal{M}}a(x)\,g(u_t)u\,
d\mathcal{M} dt.\nonumber
\end{eqnarray}

In the general case, the unique difference in the proof is that
the term $\int_0^T E(t)\,dt$ that appears on the LHS of
(\ref{4.24}) will be multiplied by a positive constant $C$,
provided that we consider $\varepsilon_i$ small enough. For
simplicity we shall assume that $C=1$.

 We will denote
\begin{eqnarray}\label{Chi}
\chi=\left[ \int_{\mathcal{M}}u_{t}\, m_{T}\cdot \nabla_T
u\,d\mathcal{M}\right]_{0}^{T}+\frac{1}{2}\left[ \int_{\mathcal{M}
}u_{t}\, u\,d\mathcal{M}\right] _{0}^{T}.
\end{eqnarray}

Next we will estimate some terms in (\ref{4.24}). Let us denote:
\begin{eqnarray}\label{R}
R:= \max_{x \in \mathcal{M}}||m(x)||_{\mathbb{R}^n}=\max_{x \in
\mathcal{M}}||x-x^0||_{\mathbb{R}^n}.
\end{eqnarray}

\smallskip

\noindent{\em Estimate for $I_1:= \int_{0}^{T}\int_{\mathcal{M}
}a(x)\,g(u_{t})(m_{T}\cdot\nabla_T u)d\mathcal{M} dt.$}

\smallskip

By Cauchy-Schwarz inequality, taking (\ref{R}) into account and
considering the inequality $ab\leq \frac{a^2}{4\eta} + \eta b^2$,
where $\eta$ is a positive number, we obtain
\begin{eqnarray}\label{5.2.14}
|I_1| \leq \frac{||a||_{L^\infty(\mathcal{M})}R^2}{\eta} \int_0^T
\int_{\mathcal{M}} a(x)|g(u_t)|^2d\mathcal{M} dt + 2\eta \int_0^T
E(t)\,dt.
\end{eqnarray}

\noindent{\em Estimate for $I_2= \frac{1}{2}\int_0^T
\int_{\mathcal{M}}a(x)\,g(u_t)u\, d\mathcal{M} dt.$}

Similarly we infer
\begin{eqnarray}\label{5.2.17}
|I_2| \leq
\frac{||a||_{L^\infty(\mathcal{M})}\lambda_1^{-1}}{16\eta}\int_0^T
\int_{\mathcal{M}} a(x)|g(u_t)|^2\,d\mathcal{M} dt +  2\eta
\int_0^T E(t)\,dt,
\end{eqnarray}
where $\lambda_1$ comes from the Poincaré inequality given in
(\ref{Poincare}).

\smallskip

Choosing $\eta =1/8$ and inserting (\ref{Chi}), (\ref{5.2.14})
 and (\ref{5.2.17}) into
(\ref{4.24}) yields
\begin{eqnarray}\label{5.2.18}
\frac12 \int_0^T E(t)\,dt &\leq& |\chi| +
C_1\int_0^T\int_{\mathcal{M}}a(x)\,(g(u_t))^2
d\mathcal{M} dt\\
&+&C_1\int_0^T\int_{\mathcal{M}_2}[ |\nabla_T u|^2 +
a(x)\,u_t^2]\,d\mathcal{M} dt\nonumber
\end{eqnarray}
where
\begin{eqnarray*}
C_1:= \max\left\{||a||_{L^\infty(\mathcal{M})}[2^{-1}\lambda_1^{-1}
+8\,R ^2], \,||B||R + |H|R,\,R\,|H|a_0^{-1}\right\},
\end{eqnarray*}
$||B||= \sup\limits_{x \in \mathcal M} \vert B_x \vert$, with $\vert
B_x \vert=\sup\limits_{\{v\in T_x \mathcal M;\vert v \vert=1\}}
\vert B_x v\vert$.

It remains to estimate the quantity
$\int_0^T\int_{\mathcal{M}_2}|\nabla_T u|^2\,d\mathcal{M} dt$ in
terms of the damping term
$\int_0^T\int_{\mathcal{M}}[a(x)\,|g(u_t)|^2 + a(x)\,|u_t|^2
]\,d\mathcal{M} dt$. For this purpose we have to built a ``cut-off''
function $\eta_\varepsilon$ on a specific neighborhood of $\mathcal
M_2$. First of all, define $\tilde \eta: \mathbb R \rightarrow
\mathbb R$ such that
\[
\begin{array}{ccc}
\tilde\eta(x) & = & \left\{
\begin{array}{clc}
1 & \mathrm{if} & x \leq 0 \\
(x-1)^2 & \mathrm{if} & x \in [1/2,1] \\
0 & \mathrm{if} & x>1
\end{array}
\right.
\end{array}
\]
and it is defined on $(0,1/2)$ in such a way that $\tilde \eta$ is a
non-decreasing function of class $C^1$. For $\varepsilon > 0$, set
$\tilde \eta_\varepsilon(x):=\tilde\eta(x / \varepsilon)$. It is
straightforward that there exist a constant $M$ which does not
depend on $\varepsilon$ such that
\[
\frac{\vert \tilde\eta_\varepsilon^\prime(x)
\vert^2}{\tilde\eta_\varepsilon(x)} \leq \frac{M}{\varepsilon^2}
\]
for every $x < \varepsilon$.

Now let $\varepsilon > 0$ such that
\[
\tilde \omega_{\varepsilon}:=\{x \in \mathcal
M;\mathrm{dist}(x,\bigcup_{i=1}^k \partial \mathcal M_{0i})<
\varepsilon\}
\]
is a tubular neighborhood of $\bigcup_{i=1}^k \partial \mathcal
M_{0i}$ and $\omega_{\varepsilon}:= \tilde \omega_{\varepsilon} \cup
\mathcal M_2$ is contained in $\mathcal M_{*}$. Define
$\eta_\varepsilon:\mathcal M \rightarrow \mathbb R$ as
\[
\begin{array}{ccc}
\eta_\varepsilon(x) & = & \left\{
\begin{array}{clc}
1 & \mathrm{if} & x \in \mathcal M_2 \\
\tilde\eta_\varepsilon(d(x,\mathcal M_2)) & \mathrm{if} & x \in
\omega_{\varepsilon} \backslash
\mathcal M_2 \\
0 & \mathrm{otherwise}. &
\end{array}
\right.
\end{array}
\]

It is straightforward that $\eta_\varepsilon$ is a function of class
$C^1$ on $\mathcal M$ due to the smoothness of $\partial M_2$ and
$\partial \omega_\varepsilon$. Notice also that
\begin{eqnarray}\label{eq:3.38'}
\frac{\vert \nabla_T \eta_\varepsilon
(x)\vert^2}{\eta_\varepsilon(x)}=\frac{\vert \tilde
\eta^\prime_\varepsilon (d(x,\mathcal M_2))\vert^2}{\tilde
\eta_\varepsilon(d(x,\mathcal M_2))} \leq \frac{M}{\varepsilon^2}
\end{eqnarray}
for every $x\in \omega_\varepsilon$. In particular, $\frac{\vert
\nabla_T \eta_\varepsilon \vert^2}{\eta_\varepsilon} \in
L^\infty(\omega_\varepsilon)$.

Taking $\xi=\eta_\varepsilon$ in the identity (\ref{identity II}) we
obtain
\begin{eqnarray} \label{5.3.1}
&&\int_0^T \int_{\omega_\varepsilon} \eta_\varepsilon |\nabla_T
u|^2 d\mathcal{M} dt \\
&=& -\left[ \int_{\omega_\varepsilon} u_t u \eta_\varepsilon
\,d\mathcal{M}\right]_0^T
+\int_0^T \int_{\omega_\varepsilon} \eta_\varepsilon |u_t|^2\,d\mathcal{M} \nonumber\\
&-& \int_0^T \int_{\omega_\varepsilon}u(\nabla_T u \cdot \nabla_T
\eta_\varepsilon) \,d\mathcal{M} dt - \int_0^T
\int_{\omega_\varepsilon} a(x)\, g(u_t)u \eta_\varepsilon\,
d\mathcal{M} dt.\nonumber
\end{eqnarray}

Next we will estimate terms on the RHS of (\ref{5.3.1}).

\smallskip
\noindent{\em Estimate for $K_1:= \int_0^T \int_{\omega_\varepsilon}
\eta_\varepsilon |u_t|^2\,d\mathcal{M} dt$}
\smallskip

From (\ref{eq:2.2}), since $\eta_\varepsilon \leq 1$ and
$\omega_{\varepsilon} \subset \mathcal{M}_{\ast}$, where the damping
lies, we deduce
\begin{eqnarray} \label{5.3.4}
K_1 \leq a_0^{-1} \int_0^T \int_{\mathcal{M}} a(x)\,
u_t^2\,d\mathcal{M}\,dt.
\end{eqnarray}

\smallskip
\noindent{\em Estimate for $K_2:= - \int_0^T
\int_{\omega_{\varepsilon}} a(x)\,g(u_t)u \eta_\varepsilon
\,d\mathcal{M} dt.$}
\smallskip

The Cauchy-Schwarz inequality, the inequality $ab \leq
\frac{1}{4\alpha}a^2 + \alpha b^2$ and (\ref{Poincare}) yield
\begin{eqnarray} \label{5.3.5}
|K_2|\leq
\frac{\lambda_1^{-1}||a||_{L^{\infty}(\mathcal{M})}}{4\alpha}\int_0^T\int_{\mathcal{M}}a(x)\,|g(u_t)|^2\,d\mathcal{M}
+ 2\alpha \int_0^T E(t)\,dt,
\end{eqnarray}
where $\alpha$ is a positive constant.

\smallskip
\noindent{\em Estimate for $K_3:=  \int_0^T
\int_{\omega_{\varepsilon}}u(\nabla_T u \cdot \nabla_T
\eta_\varepsilon) d\mathcal{M} dt.$}
\smallskip

Considering (\ref{eq:3.38'}) and applying Cauchy-Schwarz inequality,
we can write
\begin{eqnarray}\label{5.3.6}
|K_3|&\leq& \frac12 \int_0^T \left[
\int_{\omega_{\varepsilon}}\eta_\varepsilon |\nabla_T
u|^2\,d\mathcal{M} + \int_{\omega_{\varepsilon}} \frac{|\nabla_T
\eta_\varepsilon |^2}{\eta_\varepsilon}|u|^2\,d\mathcal{M}\right]dt\\
&\leq&\frac12 \int_0^T \left[
\int_{\omega_{\varepsilon}}\eta_\varepsilon |\nabla_T
u|^2\,d\mathcal{M} +
\frac{M}{\varepsilon^2}\int_{\omega_{\varepsilon}}
|u|^2\,d\mathcal{M}\right]dt. \nonumber
\end{eqnarray}

Combining (\ref{5.3.1})-(\ref{5.3.6}) we arrive to the following
inequality
\begin{eqnarray}\label{5.3.7}\qquad
\frac12 \int_0^T \int_{\omega_{\varepsilon}} \eta_\varepsilon
|\nabla_T u|^2 \,d\mathcal{M} dt &\leq& |\mathcal{Y}| +
\frac{\lambda_1^{-1}||a||_{L^{\infty}(\mathcal{M})}}{4\alpha}\int_0^T\int_{\mathcal{M}}a(x)\,|g(u_t)|^2\,d\mathcal{M}\\
&+& 2\alpha \int_0^T E(t)\,dt+ \frac{M}{2\varepsilon^2}\int_0^T
\int_{\omega_{\varepsilon}} |u|^2\,d\mathcal{M}\,dt
,\nonumber\\
&+& a_0^{-1} \int_0^T \int_{\mathcal{M}} a(x)\,
u_t^2\,d\mathcal{M}\,dt. \nonumber
\end{eqnarray}
where
\begin{eqnarray}\label{Z''}
\mathcal{Y}:=-\left[ \int_{\omega_{\varepsilon}} u_t u
\eta_\varepsilon \,d\mathcal{M}\right]_0^T.
\end{eqnarray}

Thus, combining (\ref{5.3.7}) and (\ref{5.2.18}), have in mind
that
\begin{eqnarray*}
\frac12 \int_0^T \int_{\mathcal{M}_2} |\nabla_T u|^2 \,d\mathcal{M}
dt \leq \frac12 \int_0^T \int_{\omega_{\varepsilon}}
\eta_\varepsilon |\nabla_T u|^2 \,d\mathcal{M} dt
\end{eqnarray*}
and choosing $\alpha=1/16C_1$ we deduce
\begin{eqnarray}\label{5.3.8}
&&\frac1{4}\int_0^T E(t)\,dt \leq  |\chi|  + 2 C_1|\mathcal{Y}|\qquad \\
&&+max\{C_1, 8 C_1^2 \lambda_1^{-1}||a||_{L^\infty(\mathcal M)},
2C_1 a_0^{-1}\}\int_0^T\int_{\mathcal{M}}[a(x)\,|g(u_t)|^2 +
a(x)\,|u_t|^2
]\,d\mathcal{M} dt\nonumber\\
&&+\frac{M C_1}{\varepsilon^2}\int_0^T \int_{\omega_{\varepsilon}}
\,|u|^2\,d\mathcal{M}\,dt.\nonumber
\end{eqnarray}

On the other hand, from (\ref{Chi}),  (\ref{Z''}) and (\ref{3.12})
the following estimate holds
\begin{eqnarray}\label{5.3.9}
|\chi|+  2C_2 |\mathcal{Y}|  &\leq&
C(E(0)+E(T))\\
&=& C\left[ 2\,E(T) + \int_0^T \int_{\mathcal{M}}
a(x)\,g(u_t)\,u_t\,d\mathcal{M}\right],\nonumber
\end{eqnarray}
where $C$ is a positive constant which depends also on $R$.

Then, (\ref{5.3.8}) and (\ref{5.3.9}) yield
\begin{eqnarray}\label{5.3.10}
T\,E(T)&\leq& \int_0^T E(t)\,dt\\
&\leq& C\,E(T) + C
\left[\int_0^T\int_{\mathcal{M}}[a(x)\,|g(u_t)|^2 +
a(x)\,|u_t|^2 ]\,d\mathcal{M} dt\right]\nonumber\\
&+& C\int_0^T \int_{\omega_{\varepsilon}}
\,|u|^2\,d\mathcal{M}\,dt,\nonumber
\end{eqnarray}
where $C$ is a positive constant which depends on $a_0, \lambda_1,
 R, \vert H \vert, ||B||$ and $ \frac{M}{\varepsilon^2}$.

Our aim is to estimate the last term on the RHS of (\ref{5.3.10}).
In order to do this let us consider the following lemma, where
$T_0$ is a positive constant which is sufficiently large for our
purpose.

\medskip
\begin{lemma}\label{lemma3.5}
Under the hypothesis of Theorem \ref{Theo. 3.1} and for all
$T>T_0$, there exists a positive constant $C(T_0,E(0))$ such that
if $(u,u_t)$ is the solution of (\ref{1.1}) with weak initial
data, we have
\begin{eqnarray}\label{3.53}
\int_0^T ||u(t)||_{L^2(\mathcal{M})}^2dt\leq C(T_0,E(0))
\left\{\int_0^T\int_{\mathcal{M}}\left(a(x)\,g^2(u_t))+ a(x)
u^2_t\right)d\mathcal{M}\,dt \right\}.
\end{eqnarray}
\end{lemma}

{\em Proof:} We argue by contradiction. For simplicity we shall
denote $u' := u_t$. Let us suppose that (\ref{3.53}) is not
verified and let $\{u_k(0),u'_k(0)\}$ be a sequence of initial
data where the corresponding solutions $\{u_k\}_{k\in \mathbb{N}}$
of (\ref{1.1}) with $E_k(0)$, assumed uniformly bounded in $k$,
verifies
\begin{eqnarray}\label{3.54}
\lim_{k \rightarrow +\infty}\frac{\int_0^T
||u_k(t)||_{L^2(\mathcal{M})}^2
dt}{\int_0^T\int_{\mathcal{M}}\left( a(x)\,g^2(u_k')+
a(x)\,u_k'^2\right)d\mathcal{M}\,dt}=+\infty,
\end{eqnarray}
that is
\begin{eqnarray}\label{3.55}
\lim_{k \rightarrow +\infty}\frac{\int_0^T\int_{\mathcal{M}}\left(
a(x)\,g^2(u_k')+ a(x)\,u_k'^2\right)d\mathcal{M}\,dt}{\int_0^T
||u_k(t)||_{L^2(\mathcal{M})}^2 dt}=0.
\end{eqnarray}

Since $E_k(t) \leq E_k(0)\leq L$, where $L$ is a positive
constant, we obtain a subsequence, still denoted by $\{u_k\}$ from
now on, which verifies the convergence:
\begin{eqnarray}
&& u_k \rightharpoonup u \hbox{ weakly in } H^1(\Sigma_T),\label{3.56}\\
&& u_k \rightharpoonup u \hbox{ weak star in } L^{\infty}(0,T;
V),\label{3.57}\\
&& u_k' \rightharpoonup u' \hbox{ weak star in } L^{\infty}(0,T;
L^2(\mathcal{M})).\label{3.58}
\end{eqnarray}

Employing compactness results we also deduce that
\begin{eqnarray}
u_k \rightarrow u \hbox{ strongly in }
L^2(0,T;L^2(\mathcal{M})).\label{3.60}
\end{eqnarray}

\medskip

At this point we will divide our proof into two cases, namely,
$u\ne 0$ and $u=0$.

\medskip

(i) Case (I): $u \ne 0.$

We also observe that from (\ref{3.55}) and (\ref{3.60}) we have
\begin{eqnarray}\label{3.62}
\lim_{k \rightarrow +\infty}\int_0^T\int_{\mathcal M}\left(
a(x)\,g^2(u_k')+ a(x)\,u_k'^2\right)d\mathcal{M}\,dt=0
\end{eqnarray}

Passing to the limit in the equation, when $k \rightarrow
+\infty$, we get,
\begin{equation}\label{3.63}
\left\{
\begin{aligned}
u_{tt}-\Delta_{\mathcal{M}} \,u&=& 0 \hbox{ on }~\mathcal{M} \times (0,T)\\
u_t &=& 0 \hbox{ on } \mathcal{M}_{\ast}\times(0,T),
\end{aligned}
\right.
\end{equation}
and for $u_t=v$, we obtain, in the distributional sense
\begin{equation*}
\left\{
\begin{aligned}
v_{tt}-\Delta_{\mathcal{M}} \,v&=& 0 \hbox{ on }~\mathcal{M} \times (0,T),\\
v &=& 0 \hbox{ on } \mathcal{M}_{\ast}\times (0,T).
\end{aligned}
\right.
\end{equation*}

From standard uniqueness results we conclude that $v\equiv 0$,
that is, $u_t=0$ Returning to (\ref{3.63}) we obtain the following
elliptic equation for a.e. $t\in (0,T)$ given by
\begin{equation*}
\left\{
\begin{aligned}
\Delta_{\mathcal{M}} \,u&=& 0 \hbox{ on }~\mathcal{M}\\
u_t &=& 0 \hbox{ on } \mathcal{M},
\end{aligned}
\right.
\end{equation*}
which implies that $u=0$, which is a contradiction.

\medskip
(ii) Case (II): $u = 0.$

\medskip

Defining
\begin{eqnarray}\label{3.66}
c_k := \left[ \int_0^T \int_{\mathcal{M}} |u_k|^2
d\mathcal{M}\,dt\right]^{1/2},
\end{eqnarray}
and
\begin{eqnarray}\label{3.67}
\overline{u}_k:= \frac{1}{c_k}\,u_k,
\end{eqnarray}
we obtain
\begin{eqnarray}\label{3.68}
&\int_0^T\int_{\mathcal{M}} |\overline{u}_k|^2 d\mathcal{M}\,dt =
\int_0^T\int_{\mathcal{M}} \frac{|u_k|^2}{c_k^2}d\mathcal{M}\,dt=
\frac{1}{c_k^2}\int_0^T \int_{\mathcal{M}} |u_k|^2
d\mathcal{M}\,dt=1. &
\end{eqnarray}

Setting
\begin{eqnarray*}
\overline{E}_k(t)&:=& \frac12 \int_{\mathcal{M}}
|\overline{u}'_k|^2\,d\mathcal{M} + \frac12 \int_{\mathcal{M}}
|\nabla \overline{u}_k|^2\,d\mathcal{M} ,
\end{eqnarray*}
we deduce automatically that
\begin{eqnarray}\label{3.69}
\overline{E}_k(t)= \frac{E_k(t)}{c_k^2}.
\end{eqnarray}

Recalling (\ref{5.3.10}) we obtain, for  $T$ large enough that
\begin{eqnarray*}
E(T) \leq \hat{C}\left[\int_0^T \int_{\mathcal{M}} (a(x)\,
g^2(u_t)+ a(x)\,u_t^2 )\,d\mathcal{M}\,dt +
\int_0^T\int_{\mathcal{M}} |u|^2\,d\mathcal{M}\,dt\right],
\end{eqnarray*}
and employing the identity $E(T)-E(0)= -\int_0^T\int_{\mathcal{M}}
a(x)\,g(u_t)\,u_t\,d\mathcal{M}\,dt$, we can write
\begin{eqnarray*}
E(t) \leq E(0) \leq \tilde{C}\left[\int_0^T \int_{\mathcal{M}}
(a(x)\,g^2(u_t)+ a(x)\,u_t^2 )\,d\mathcal{M}\,dt +
\int_0^T\int_{\mathcal{M}} |u|^2\,d\mathcal{M}\,dt\right],
\end{eqnarray*}
for all $t\in (0,T)$, with $T$ large enough. The last inequality
and (\ref{3.69}) give us
\begin{eqnarray}\label{3.70}
\overline{E}_k(t):=\frac{E_k(t)}{c_k^2} \leq \tilde{C}
\left[\frac{\int_0^T \int_{\mathcal{M}}(a(x)\,g^2(u_k')+
a(x)\,u_k'^2)}{\int_0^T
\int_{\mathcal{M}}|u_k|^2\,d\mathcal{M}\,dt}+1 \right].
\end{eqnarray}

From (\ref{3.55}) and (\ref{3.70}) we conclude that there exists a
positive constant $\hat{M}$ such that
\begin{eqnarray*}
\overline{E}_k(t):=\frac{E_k(t)}{c_k^2} \leq \hat{M}, ~\hbox{ for
all } t\in[0,T] ~\hbox{ and for all } k\in \mathbb{N},
\end{eqnarray*}
that is,
\begin{eqnarray}\label{3.71}
\frac12 \int_{\mathcal{M}} |\overline{u}'_k|^2\,d\mathcal{M} +
\frac12 \int_{\Omega} |\nabla \overline{u}_k|^2\,d\mathcal{M}\leq
\hat{M}, ~\hbox{ for all } t\in[0,T] ~\hbox{ and for all } k\in
\mathbb{N}.
\end{eqnarray}

For a subsequence $\{\overline{u}_k\}$, we obtain
\begin{eqnarray}
&&\overline{u}_k \rightharpoonup \overline{u} \hbox{ weak star
in }L^{\infty}(0,T;V),\label{3.72}\\
&&\overline{u}'_k  \rightharpoonup \overline{u}' \hbox{ weak
star in } L^{\infty}(0,T;L^2(\mathcal{M})),\label{3.73}\\
&&\overline{u}_k  \rightarrow \overline{u} \hbox{ strongly in }
L^2(0,T;L^2(\mathcal{M})).\label{3.74}
\end{eqnarray}

We observe that from (\ref{3.55}) we deduce
\begin{eqnarray}\label{3.75}
\lim_{k \rightarrow +\infty}\int_0^T \int_{\mathcal{M}}
\frac{a(x)\,g^2(u_k')}{c_k^2}\,d\mathcal{M}\,dt=0~\hbox{ and }
\lim_{k \rightarrow +\infty}\int_0^T \int_{\mathcal{M}}
a(x)\,|\overline{u}_k'|^2\,d\mathcal{M}\,dt=0.
\end{eqnarray}

In addition $\overline{u}_k$ satisfies the equation
\begin{equation*}
\left.
\begin{aligned}
\overline{u}_k'' - \Delta_{\mathcal{M}}\overline{u}_k +
a(x)\,\frac{g(\overline{u}_k')}{c_k}&=&0\quad \hbox{ on }
\mathcal{M} \times (0,T).
\end{aligned}
\right.
\end{equation*}

Passing to the limit when $k \rightarrow +\infty$ taking the above
convergence into account, we obtain
\begin{equation}\label{P2}
\left\{
\begin{aligned}
\overline{u}'' - \Delta_{\mathcal{M}}\overline{u} &=&0\quad
\hbox{ on } \mathcal{M} \times (0,T),\\
\overline{u}'&=&0\quad \hbox{ on }\mathcal{M}_{\ast} \times (0,T).
\end{aligned}
\right.
\end{equation}

Then, $v= \overline{u}_t$ verifies, in the distributional sense
\begin{equation*}
\left\{
\begin{aligned}
v_{tt}-\Delta_{\mathcal{M}} \,v&=& 0 \hbox{ on }~\mathcal{M}\\
v &=& 0 \hbox{ on } \mathcal{M}_{\ast}.
\end{aligned}
\right.
\end{equation*}

Applying uniqueness standard results it results that
$v=\overline{u}_t=0$. Returning to (\ref{P2}) we obtain, for a.e.
$t\in (0,T)$ that
\begin{equation*}
\left\{
\begin{aligned}
\Delta_{\mathcal{M}} \,\overline{u}&=& 0 \hbox{ on }~\mathcal{M}\\
\overline{u}_t &=& 0 \hbox{ on } \mathcal{M},
\end{aligned}
\right.
\end{equation*}
from what we deduce that $\overline{u}=0$, which is a
contradiction in view of (\ref{3.68}) and (\ref{3.74}). The lemma
is proved. \quad $\square$

\smallskip

Inequalities (\ref{5.3.10}) and (\ref{3.53}) lead us to the
following result.
\smallskip

\noindent \textbf{Proposition 5.2.2:} \textit{For }$T>0$\textit{\
large enough, the solution }$\left[ u,u_{t}\right] $\textit{\ of
}(\ref{3.1})
\textit{satisfies}%
\begin{equation}
E(T)\leq C\,\int_{0}^{T}\int_{\mathcal{M} _{2}}\left[ a(x)\left\vert
u_{t}\right\vert ^{2}+ a(x)\left\vert g\left( u_{t}\right)
\right\vert ^{2}\right] d\mathcal{M} dt \label{5.3.18}
\end{equation}%
\textit{where the constant }$C=C(T_0, E(0), a_0, \lambda_1, R,
||B||, \frac{M}{\varepsilon^2}).$

\subsection{Conclusion of Theorem 3.1}

In what follows we will proceed exactly as in Lasiecka and
Tataru's  work\cite{Lasiecka-Tataru}(see Lemma 3.2 and Lemma 3.3
of the referred paper) adapted o our context. Let
\begin{eqnarray*}
\Sigma _{\alpha } &=&\left\{ \left( t,x\right) \in \Sigma _{1}/\text{ \ }%
\left\vert u_{t}\right\vert >1\text{ \ a. e.}\right\} , \\
\Sigma _{\beta } &=&\Sigma _{1}\backslash \Sigma _{\alpha }.
\end{eqnarray*}

Then using hypothesis $(iii)$ in Assumption \ref{ass:1}, we obtain%
\begin{equation}
\int_{\Sigma _{\alpha }}a(x)\left( \left[ g\left( u_{t}\right)
\right] ^{2}+\left( u_{t}\right) ^{2}\right) d\Sigma _{\alpha
}\leq \left( k^{-1}+K\right) \int_{\Sigma _{\alpha }}a(x)g\left(
u_{t}\right) u_{t}d\Sigma _{\alpha }.  \label{5.4.1}
\end{equation}%

Moreover, from (\ref{4.1'})
\begin{equation}
\int_{\Sigma _{\beta }}a(x)\left( \left[ g\left( u_{t}\right)
\right] ^{2}+\left( u_{t}\right) ^{2}\right) d\Sigma _{\beta }\leq
(1+||a||_{\infty})\int_{\Sigma _{\beta }}h\left(a(x) g\left(
u_{t}\right) u_{t}\right) d\Sigma _{\beta }. \label{5.4.2}
\end{equation}%

Then by Jensen's inequality
\begin{eqnarray}
(1+||a||_{\infty})\int_{\Sigma _{\beta }}h\left( g\left(
u_{t}\right) u_{t}\right) d\Sigma _{\beta } &\leq
&(1+||a||_{\infty})meas\left( \Sigma\right) h\left(
\frac{1}{meas\left( \Sigma\right) }\int_{\Sigma}a(x)g\left(
u_{t}\right) u_{t}d\Sigma
\right) \smallskip  \notag \\
&=&(1+||a||_{\infty})meas\left( \Sigma\right) r\left(
\int_{\Sigma}a(x)g\left( u_{t}\right) u_{t}d\Sigma\right) ,
\label{5.4.3}
\end{eqnarray}%
where $r\left( s\right) =h\left( \frac{s}{meas\left( \Sigma \right) }%
\right)$ is defined in (\ref{4.2'}). Thus%
\begin{eqnarray}
\int_{\Sigma}a(x)\left( \left[ g\left( u_{t}\right) \right]
^{2}+\left( u_{t}\right) ^{2}\right) d\Sigma &\leq &\left(
k^{-1}+K\right)
\int_{\Sigma}a(x) g\left( u_{t}\right) _{t}d\Sigma _{1}\smallskip  \notag \\
&&+(1+||a||_{\infty})meas\left( \Sigma\right) r\left( \int_{\Sigma
_{1}}a(x)g\left( u_{t}\right) u_{t}d\Sigma\right) . \label{5.4.4}
\end{eqnarray}

Splicing, together, (\ref{5.3.18}) and (\ref{5.4.4}), we have%
\begin{eqnarray}
E(T) &\leq &(1+||a||_{\infty})C\,\left[
\frac{K_{0}}{(1+||a||_{\infty})}\int_{\Sigma}g\left( u_{t}\right)
u_{t}d\Sigma _{1}\right. \smallskip  \notag \\
&&\left. + meas\left( \Sigma\right) r\left(
\int_{\Sigma}a(x)\,g\left( u_{t}\right) u_{t}d\Sigma\right)
\right] , \label{5.4.5}
\end{eqnarray}
where $K_{0}=k^{-1}+K$. Setting%
\begin{eqnarray*}
L &=&\frac{1}{C\,meas\left( \Sigma\right) (1+||a||_{\infty})},\smallskip \\
c &=&\frac{M_{0}}{meas\left( \Sigma \right)(1+||a||_{\infty}) },
\end{eqnarray*}%
we obtain
\begin{eqnarray}
p\left[ E(T)\right] &\leq &\int_{\Sigma}a(x)\,g\left( u_{t}\right)
u_{t}\,d\Sigma  =E(0)-E(T),  \label{5.4.6}
\end{eqnarray}%
where the function $p$ is as defined in (\ref{4.3'}). To finish
the proof of Theorem 3.1, we invoke the following result from I.
Lasiecka et al. \cite{Lasiecka-Tataru}:

\smallskip

\noindent \textbf{Lemma B}:\textit{\ Let }$p$\textit{\ be a positive,
increasing function such that }$p(0)=0$\textit{. Since }$p$\textit{\ is
increasing we can define an increasing function }$q,$ $q(x)=x-(I+p)^{-1}%
\left( x\right) .$\textit{\ Consider a sequence }$s_{n}$\textit{\ of
positive numbers which satisfies}%
\begin{equation*}
s_{m+1}+p(s_{m+1})\leq s_{m}.
\end{equation*}%
\textit{Then }$s_{m}\leq S(m)$\textit{, where }$S(t)$\textit{\ is a solution
of the differential equation}%
\begin{equation*}
\frac{d}{dt}S(t)+q(S(t))=0,\text{ \ }S(0)=s_{0}.
\end{equation*}%
\textit{Moreover, if }$p(x)>0$\textit{\ for }$x>0$\textit{, then }$\underset{%
t\rightarrow \infty }{\lim }$\textit{\ }$S(t)=0.$

With this result in mind, we replace $T$ (resp. $0$) in (\ref{5.4.6}) with $%
m(T+1)$ (resp. \ $mT$) to obtain%
\begin{equation}
E(m(T+1))+p\left( E(m(T+1))\right) \leq E(mT),\text{ \ for \
}m=0,1,.... \label{5.4.7}
\end{equation}

Applying Lemma B with $s_{m}=E(mT)$ thus results in%
\begin{equation}
E(mT)\leq S(m),\text{ \ \ }m=0,1,....  \label{5.4.8}
\end{equation}%

Finally, using the dissipativity of $E(t)$ inherent in the
relation (\ref{3.12})
, we have for $t=mT+\tau ,$ $0\leq \tau \leq T,$%
\begin{equation}
E(t)\leq E(mT)\leq S(m)\leq S\left( \frac{t-\tau }{T}\right) \leq
S\left( \frac{t}{T}-1\right) \text{ \ \ for \ }t>T\text{,}
\label{5.4.9}
\end{equation}%
where we have used above the fact that $S(.)$ is dissipative. The
proof of Theorem 3.1 is now completed.

\subsection{Further Remarks}

From the above procedure, we can construct a wide assortment of
compact surfaces by jointing pieces of different kind of surfaces.
However, according to the position of the observer (point $x_0$)
the dissipative and the non dissipative areas can change
drastically. To illustrate this, let us consider the Catenoid or
the Trinoid (see figures 2 and 3 below) that are minimal surfaces,
that is $H=0$.

{\small
\begin{figure}[ht]
\begin{picture}(8,6.5)(-0.5,-3)

 \qbezier(1,-1)(1.2,-2)(3.35,-2)   \qbezier(3.35,-2)(5.5,-2)(5.7,-1)
\qbezier(1,-1)(1,0)(2.2,0.1)  \qbezier(4.5,0.1)(5.7,0)(5.7,-1)
\qbezier(1,-1)(3.5,0)(1,1)  \qbezier(5.7,-1)(3.2,0)(5.7,1)
\qbezier(1,1)(1.5,0)(3.35,0) \qbezier(3.35,0)(5.2,0)(5.7,1)
\qbezier(1,1)(1.5,2)(3.35,2) \qbezier(3.35,2)(5.2,2)(5.7,1)
\put(3.35,2){\line(0,-1){1.5}}\put(3.30,0){\line(0,-1){2}}\put(2.85,0.25){$H=0$}
\qbezier(2.3,1.9)(2.9,1.3)(2.8,0.05)
\qbezier(1.5,1.6)(2.3,1.3)(2.45,0.1)
\qbezier(4.40,1.9)(3.9,1.3)(3.8,0.05)
\qbezier(5.4,1.4)(4.3,1.3)(4.15,0.05)
\qbezier(2.7,0)(3,-1)(2.4,-1.9)\qbezier(2.4,0.05)(2.5,-1)(1.5,-1.6)
\qbezier(3.85,0)(3.6,-1)(4.1,-1.95)\qbezier(4.3,0.1)(4.1,-1)(5.4,-1.5)
\qbezier(2,0.2)(1.2,1)(2,1.3)\qbezier(5,1.3)(5.5,1)(4.7,0.2)\qbezier(2,1.3)(3.35,2)(5,1.3)
\qbezier(2.7,0.08)(1.7,0.8)(2.5,1.1)
\qbezier(4.3,1.1)(5,0.8)(4.1,0.08)
\qbezier(2.5,1.1)(3.4,1.5)(4.3,1.1)
\qbezier(3.2,0.03)(2.2,0.6)(3,0.8)
\qbezier(3.5,0.03)(4.5,0.6)(4,0.8) \qbezier(3,0.8)(3.35,1)(4,0.8)
\qbezier(2.2,-0.1)(1.2,-0.5)(1.7,-1)\qbezier(4.5,-0.1)(5.2,-0.5)(4.7,-1)\qbezier(1.7,-1)(3.2,-1.8)(4.7,-1)
\qbezier(2,-0.5)(3.3,-1.1)(4.65,-0.4)



\put(2,-3){Figure 2: Catenoid}

\end{picture}
\bigskip
\end{figure}}

{\small
\begin{figure}[ht]
\begin{picture}(8,6.5)(-0.5,-3)

\qbezier(2.7,-3)(2.5,-1.5)(3.5,-1.5)
\qbezier(3.5,-1.5)(4.5,-1.5)(4.3,-3)
\qbezier(4.8,-0.5)(4.8,0.5)(6,0.7)
\qbezier(4.8,-0.5)(5.3,-1.5)(6.7,-0.5)
\qbezier(1.8,-0.5)(1.8,0.5)(0.6,0.7)
\qbezier(-0.1,-0.5)(1.3,-1.5)(1.8,-0.5)
\qbezier(0.6,0.7)(2.7,0.9)(6,0.7)
\qbezier(4.3,-3)(5.1,-1.4)(6.7,-0.5)
\qbezier(4.3,-3)(5.4,-2.5)(6.7,-0.5)
\qbezier(-0.1,-0.5)(1.6,-1.5)(2.7,-3)
\qbezier(-0.1,-0.5)(1.2,-2.2)(2.7,-3)
\put(3,-0.5){$H=0$} \put(2,-4){Figure 4: Trinoid}
\end{picture}
\bigskip
\end{figure}}

{\small
\begin{figure}[ht]
\begin{picture}(8,6.5)(-0.5,-3)

\qbezier(-2,-1)(0,-2)(2,-1) \qbezier(-0.5,0.8)(0,0.6)(0.5,0.8)
\qbezier(-2,-1)(-0.8,0)(-0.5,0.8)\qbezier(0.5,0.8)(0.8,0)(2,-1)
\qbezier(-0.5,0.8)(-0.4,1)(-0.5,1.2)\qbezier(0.5,0.8)(0.4,1)(0.5,1.2)
\qbezier(-0.5,1.2)(-1,1.8)(0,2)\qbezier(0,2)(1,1.8)(0.5,1.2)

\qbezier(4,1)(6,2)(8,1)\qbezier(5.5,-0.8)(6,-0.5)(6.5,-0.8)
\qbezier(4,1)(5.2,0)(5.5,-0.8)\qbezier(6.5,-0.8)(6.8,0)(8,1)
\qbezier(5.5,-0.8)(5.6,-1)(5.5,-1.2)\qbezier(6.5,-0.8)(6.4,-1)(6.5,-1.2)
\qbezier(5.5,-1.2)(5,-1.8)(6,-2)\qbezier(6.5,-1.2)(7,-1.8)(6,-2)
\qbezier(-2,-1)(-2.5,-1.5)(-2,-2)\qbezier(2,-1)(2.5,-1.5)(2,-2)\qbezier(-2,-2)(0,-3.5)(2,-2)
\qbezier(4,1)(3.5,1.5)(4,2)\qbezier(8,1)(8.5,1.5)(8,2)\qbezier(4,2)(6,3.5)(8,2)
\put(-0.5,-0.5){$H=0$}\put(-0.5,-2.2){$H\leq
0$}\put(-0.5,1.4){$H\leq
0$}\put(0,-3.5){$x_0$}\put(0,-3.2){$\centerdot$}\put(6,-3){$x_0$}\put(6,-2.7){$\centerdot$}
\put(-.5,2.5){FIG. $A$}\put(3.5,2.5){FIG. $B$}
\put(5.5,0.5){$H=0$}\put(5.5,-1.5){$H\leq0$}\put(5.5,2.0){$H\leq0$}

\put(2,-3.8){Figure 3.}
\end{picture}
\bigskip
\end{figure}}

\bigskip
Considering a strategic piece of one of these surfaces we can
construct another compact surface  according to the figure 5 above.
Remember that the {\em non dissipative regions} must occur where
$m(x)\cdot \nu(x) \leq 0$, $H \leq 0$ and simultaneously we are
forced to consider $|k_1-k_2|$ sufficiently small (by parts). The
{\em dissipative area} must contain strictly the closure of the
points $x\in \mathcal{M}$ such that $m(x) \cdot \nu(x)>0$. It is not
difficult to see that the non dissipative area in the figure $A$ can
be located near the top and/or near $x_0$ while the non dissipative
area in the figure $B$ can be located in the middle in the middle of
surface and/or near $x_0$, assuming evidently that $k_1\approx k_2$
on these non dissipative areas.

\end{document}